\documentclass{custom}
\usepackage{leftidx}
\usepackage{quiver}
\usepackage{amsthm}
\usepackage{thmtools}
\usepackage[capitalize]{cleveref}
\usepackage[utf8]{inputenc}
\usepackage[greek,english]{babel}
\usepackage{tikz}
\usepackage{etoolbox}
\makeatletter
\patchcmd{\@maketitle}
  {\ifx\@empty\@dedicatory}
  {\ifx\@empty\@date \else {\vskip3ex \centering\footnotesize\@date\par\vskip1ex}\fi
   \ifx\@empty\@dedicatory}
  {}{}
\patchcmd{\@adminfootnotes}
  {\ifx\@empty\@date\else \@footnotetext{\@setdate}\fi}
  {}{}{}
\makeatother

\citesort

\DeclareMathOperator{\End}{End}
\DeclareMathOperator{\GL}{GL}
\DeclareMathOperator{\Ann}{Ann}
\DeclareMathOperator{\Jac}{Jac}
\DeclareMathOperator{\SU}{SU}
\newcommand{\lMod}[1]{\mathsf{\leftidx_{\mathnormal{#1}}Mod}}
\newcommand{\lAnn}[1]{\opn{\leftidx_{\mathnormal{#1}}{Ann}}}
\newcommand{\lProj}[1]{\mathsf{\leftidx_{\mathnormal{#1}}Proj}}
\newcommand{\Vect}{\mathsf{Vect}}
\newcommand{\K}{\mathbb{K}}
\newcommand{\Z}{\mathbb{Z}}
\newcommand{\h}{\mathcal{H}}
\newcommand{\A}{\mathcal{A}}
\newcommand{\D}{\mathcal{D}}
\newcommand{\?}{{?}}

\newcommand{\opn}[1]{\operatorname{#1}}

\newtheorem{thm}{Theorem}[section]
\newtheorem{prop}[thm]{Proposition}
\newtheorem{cor}[thm]{Corollary}
\newtheorem{lemma}[thm]{Lemma}
\theoremstyle{definition}
\newtheorem{defn}[thm]{Definition}
\newtheorem{rem}[thm]{Remark}
\newtheorem{example}[thm]{Example}
\Crefname{lemma}{Lemma}{Lemmata}

\numberwithin{equation}{section}

\crefformat{section}{\S#2#1#3}
\crefmultiformat{section}{\S\S#2#1#3}{ and~#2#1#3}{, #2#1#3}{, and~#2#1#3}

\title{McKay quivers of Hopf algebras}

\author{James Jensen}
\date{November 2026}

\begin{document}

\begin{abstract}
We present results concerning the connectivity of McKay quivers of finite-dimensional Hopf algebras with no assumption of semisimplicity.
Applying these results to group algebras, we give the number of connected components in the McKay quiver of a group over a field of prime characteristic.
\end{abstract}

\maketitle

\section{Introduction}

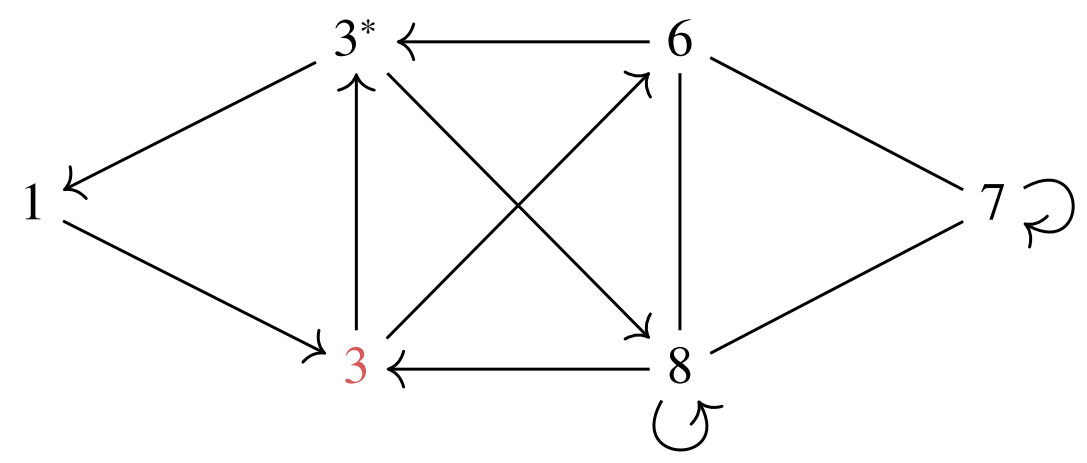
\begin{figure}
	\[\begin{tikzcd}[cramped]
		&& {3^*} && 6 && \\
		1 &&&&&& 7 \\
		&& \textcolor{rgb,255:red,214;green,92;blue,92}{3} && 8
		\arrow[from=1-3, to=2-1]
		\arrow[from=1-3, to=3-5]
		\arrow[from=1-5, to=1-3]
		\arrow[no head, from=1-5, to=2-7]
		\arrow[from=2-1, to=3-3]
		\arrow[from=2-7, to=2-7, loop, in=330, out=30, distance=5mm]
		\arrow[no head, from=2-7, to=3-5]
		\arrow[from=3-3, to=1-3]
		\arrow[from=3-3, to=1-5]
		\arrow[no head, from=3-5, to=1-5]
		\arrow[from=3-5, to=3-3]
		\arrow[from=3-5, to=3-5, loop, in=300, out=240, distance=5mm]
	\end{tikzcd}\]
	\caption{
		The McKay quiver of $\opn{PSL}_3(2)$ at the module represented by the red vertex. 
		Vertices are labeled by their dimension.
		We depict a pair of opposite directed edges as a single undirected edge.
		Note that the dimensions form a left and right eigenvector of eigenvalue 3.
	}
	\label{fig:PSL32}
\end{figure}

Given a group $G$ and $\mathbb C G$-module $V$, the \emph{McKay quiver} $\Gamma_{\mathbb CG}(V)$ of $G$ at $V$ is a directed graph describing how tensor products of $\mathbb C G$-modules with $V$ decompose.
The vertices of $\Gamma_{\mathbb CG}(V)$ are the (isomorphism classes of the) simple $\mathbb CG$-modules $S_1, \dots, S_t$,
and we have an edge from $S_i$ to $S_j$ for every copy of $S_j$ appearing as a direct summand of $S_j \otimes V$.
Therefore if $(M_{ij})$ is the adjacency matrix of the quiver,
called the \emph{McKay matrix} of $G$ at $V$,
we have
\begin{equation*}
	S_i \otimes V \cong_{\mathbb C G} \bigoplus_{j=1}^t S_j^{\oplus M_{ij}}.
\end{equation*}
This construction, an example of which is shown in \cref{fig:PSL32}, was first proposed by John McKay in \cite{mckayGraphsSingularitiesFinite1980},
where it was observed that if $G$ is a finite subgroup of $\SU(2)$ and $V = \mathbb{C}^2$ equipped with the natural action of $G$,
then the corresponding quiver is one of the affine Dynkin diagrams of type $\tilde{A}_n, \tilde{D}_n$, or $\tilde{E}_n$.
This empirical observation sets the stage for a bijection between conjugacy classes of finite subgroups of $\SU(2)$ on the one hand, and simply-laced Dynkin diagrams on the other, known as the \emph{McKay correspondence}.
Connections to other areas of mathematics abound,
as the same Dynkin diagrams classify regular polyhedra, certain semisimple Lie algebras, and Kleinian singularities, among others~\cite{buchweitzPlatonicSolidsPreprojective2013}.

A fundamental result about McKay quivers, already stated in \cite{mckayGraphsSingularitiesFinite1980},
is that the McKay quiver of $G$ at $V$ is strongly connected whenever $V$ is faithful on $G$,
a consequence of Burnside's theorem~{\cite[\S XV, Theorem~IV]{burnsideTheoryGroupsFinite2012}}:

\begin{thm} \label{thm:burnside-classical}
	If $G$ is a finite group and $V$ be a $\mathbb C G$-module faithful on $G$,
	then every simple $\mathbb CG$-module appears as a direct summand of $V^{\otimes n}$ for some $n \in \mathbb N$.
\end{thm}

More recently, the connectivity properties of McKay quivers were expanded upon by Browne~\cite{browneConnectivityPropertiesMcKay2021},
who showed that the strongly and weakly connected components of McKay quivers coincide,
and characterised the connected components of McKay quivers at non-faithful $\mathbb C G$-modules.

The setup for McKay quivers makes sense anywhere one can meaningfully talk about tensor products and their decomposition into ``atoms.''\footnote{The natural categorical environment for this is that of \emph{fusion} or more generally \emph{tensor categories} in the sense of \cite{etingofTensorCategories2015}, although we do not consider McKay quivers in this generality here.}
Let $\K$ denote a field.
For a unital associative $\K$-algebra $\h$, hereinafter called a \emph{$\K$-monoid},
a well-behaved tensor product on $\h$-modules is equivalent to $\h$ bearing the additional structure of a $\K$-bimonoid.
A special subclass of these are Hopf $\K$-algebras,
whose representations admit duals akin to contragredient representations of groups.
These definitions are reviewed in \cref{sec:prelim}.

The study of McKay quivers in this generality is fairly recent.
A series of papers by Chan et al.~\cite{chanQuantumBinaryPolyhedral2016,chanMcKayCorrespondenceSemisimple2018,chanMcKayCorrespondenceSemisimple2019}  proved a quantum anologue of the classical McKay correspondence for Hopf algebras called quantum binary polyhedral groups.
McKay quivers of general Hopf $\K$-algebras were considered by Grinberg et al.~\cite{grinbergCriticalGroupsHopf2020} in the context of the chip firing game,
where it was observed that the dimensions of simple modules give a right eigenvector to the McKay matrix.
In a similar vein, a paper by Benkart et al.~\cite{benkartMcKayMatricesFinitedimensional2022} explored the statistical properties of certain Markov chains constructed from McKay quivers, including those of Hopf $\K$-algebras.
A later paper by Benkart et al.~\cite{benkartMcKayMatricesFinitedimensional2022} gives a more thorough treatment of the spectral properties of McKay matrices of finite-dimensional Hopf $\K$-algebras,.
Finally, it was shown by Crawford~\cite{crawfordSuperpotentialsQuiverAlgebras2023}
that for a semisimple finite-dimensional Hopf $\K$-algebra $\h$, the McKay quiver of $\h$ at $V$ is strongly connected if and only if $V$ is \emph{inner faithful},
generalising McKay's connectivity criterion.

In this paper we present connectivity results for the McKay quivers of possibly non-semisimple Hopf $\K$-algebras.
In particular we show that strongly and weakly connected components coincide (\cref{thm:mckay-strong-weak}) and give a sharp criterion for connectivity (\cref{thm:criterion}).
Since these results do not require semisimplicity,
they may be applied to McKay quivers of groups over modular fields.
\Cref{thm:modular-mckay-components} generalises a result of Browne \cite{browneConnectivityPropertiesMcKay2021} giving the number of connected components to this setting.

\subsection*{Conventions}
The set of nonnegative integers is denoted $\mathbb N$, and $\K$ denotes a field.
Unadorned tensor products are over $\K$,
and dimension is $\K$-dimension.
All $\K$-spaces, including $\mathbb K$-monoids and their modules, are finite-dimensional,
and all modules are left modules.
For a $\K$-monoid $\A$,
we denote its category of modules $\lMod \A$,
and the symbols $\cong_\A$ and $\leq_\A$ then represent the isomorphism and submodule relations of $\A$-modules.
An $\A$-morphism is a homomorphism of $\A$-modules.          

\section{Preliminaries\label{sec:prelim}}

Proofs concerning Grothendieck groups can be found in \cite{lorenzRepresentationsFinitedimensionalHopf1997}.
Let $\A$ be a $\K$-monoid
with multiplication $\mu : \A \otimes \A \to \A$
and unit $\eta : \K \to \A$.
Let $S_1, \dots, S_t$ be the non-isomorphic simple $\A$-modules and $d_i = \dim \End_\A(S_i)$,
so that $\K$ is a splitting field for $\A$ if and only if all $d_i = 1$.
The projective covers $P_1,\dots,P_t$ of the non-isomorphic simple modules are precisely the indecomposable projective modules, i.e. the indecomposable direct summands of the free module $\A$, where we have
\begin{align} \label{eq:regular remak}
	P_i / \opn{rad}P_i &\cong S_i,
	&
	\A &\cong_{\A} \bigoplus_{i=1}^t P_i^{\oplus \dim S_i}.
\end{align}
From the abelian category $\lMod \A$ and its full additive subcategory $\lProj \A$ spanned by projective modules, we can construct two \emph{Grothendieck groups}, 
algebraic gadgets which track the Jordan-Hölder and Remak decompositions of modules and projective modules respectively.

\begin{enumerate}
	\item $\mathrm G_0\A$, the Grothendieck group of $\lMod \A$,
		is the free abelian group generated by isomorphism classes of $\A$-modules
		subject to the relation $[U] - [V] + [W] = 0$ for every exact sequence $0 \to U \hookrightarrow V \twoheadrightarrow W \to 0$.
		The simple modules $[S_i]$ form a $\Z$-basis so that
		\begin{equation*}
			[\A] = \sum_{j=1}^t (\dim P_j) [S_j] \in \mathrm G_0 \A.	
		\end{equation*}
		We denote the coefficient of $[S_i]$ in the basis expansion of an $\A$-module $V$ as $[V : S_i]$.
		Note that $[V : S_i] \geq 1$ if and only if $S_i$ is a subquotient of $V$.
	\item $\mathrm K_0 \A$, the Grothedieck group of $\lProj \A$,
		is the free group generated by isomorphism classes of projective $\A$-modules subject to the relation $[U] - [V] + [W] = 0$ for every exact sequence $0 \to U \hookrightarrow V \twoheadrightarrow W \to 0$,
		or what is the same by projectivity, whenever $V \cong_\A U \oplus W$.
		The indecomposable projective modules $[P_i]$ form a $\Z$-basis so that
		\begin{equation*}
			[\A] = \sum_{i=1}^t (\dim S_i) [P_i] \in \mathrm K_0 \A.	
		\end{equation*}
		We denote the coefficient of $[P_i]$ in the basis expansion of a projective $\A$-module $P$ as $[P : P_i]$.
\end{enumerate}
When there is no possibility of confusion, we abuse notation and write $V$ for the corresponding element $[V]$ of $\mathrm G_0 \A$ or $\mathrm K_0 \A$.

The inclusion $\lProj \A \hookrightarrow \lMod \A$ gives a group homomorphism $c : \mathrm K_0 \A \to \mathrm G_0\A$ called the \emph{Cartan map},
and the dimension of $\A$-modules gives well-defined group homomorphisms $\mathrm G_0 \A \to \Z$.
There is also well-defined bilinear pairing
\begin{align*}
	\langle \? , \? \rangle : \mathrm K_0 \A \times \mathrm G_0 \A &\to \Z \\
	(P, V) &\mapsto  \dim \lMod \A(P, V).
\end{align*}
where in particular $\langle P_i, S_j\rangle = \delta_{ij}d_j$ and thus $\langle P_i, V \rangle = [V : S_i]d_i$ for any module $V$.

In order to take tensor products of $\A$-modules, we require additional structure.

\begin{defn}[{See \cite[\S 3]{sweedlerHopfAlgebras1969}}]
	A \emph{$\K$-bimonoid} $\h$ is a $\K$-monoid equipped with $\K$-monoid homomorphisms
	\begin{align*}
		\Delta : \h &\to \h \otimes \h,
		&
		\epsilon : \h \to \K
	\end{align*}
	termed the \emph{comultiplication} and \emph{counit} respectively,
	where these data satisfy
	\begin{align*}
		(\Delta \otimes 1) \circ \Delta &= (1 \otimes \Delta) \circ \Delta,
		&
		(\epsilon \otimes 1) \circ \Delta &= 1 = (1 \otimes \epsilon) \circ \Delta
	\end{align*}
	called \emph{coassociativity} and \emph{counitality} respectively.\footnote{
		We use $(=)$ for equality up to unitor and associator morphisms in $\Vect_\K$.
	}
	A $\K$-bimonoid is called \emph{cocommutative} if and only if $\tau \circ \Delta = \Delta$, where $\tau : x \otimes y \mapsto y \otimes x$ is the canonical braiding of $\lMod \K$.
\end{defn}

\begin{example}\label{ex:monoid}
	Given a finite monoid $M$, we can construct two $\K$-bimonoids.
	The \emph{monoid algebra} $\K M$ is a cocommutative $\K$-bimonoid with $M$ as a basis and
	\begin{align*}
		\mu &: x \otimes y \mapsto xy ,
		&
		\eta &: 1 \mapsto 1_M,
		\\
		\Delta &: x \mapsto x \otimes x,
		&
		\epsilon &: x \mapsto 1
	\end{align*}
	for $x, y \in M$.
	The \emph{algebra of functions} $\K^M$ is a commutative $\K$-bimonoid whose elements are arbitrary maps $M \to \K$, with
	\begin{align*}
		\mu &: f \otimes g \mapsto x \mapsto f(x)g(x)	
		&
		\eta &: 1 \mapsto m \mapsto 1
		\\
		\Delta &: f \mapsto (x, y) \mapsto f(y)
		&
		\epsilon &: f \mapsto f(1)
	\end{align*}
	where we have abused notation to make the canonical isomorphism $\K^M \otimes \K^M \cong_\K \K^{M \times M}$ implicit.
\end{example}

We employ a notation due to Sweedler~\cite{sweedlerHopfAlgebras1969} and write
\begin{equation*}
	\Delta(x) =: \sum_{(x)} x_{(1)} \otimes x_{(2)},
\end{equation*}
a nod to the fact that $\Delta(x)$ can be written as a sum of simple tensors
while avoiding choosing such a decomposition.\footnote{The $x_{(i)}$s may then only be manipulated such that the result is independent of such a choice, i.e. multilinearly.}
Coassociativity means we may unambiguously write
\begin{align*}
	\Delta_n : \h &\to \h^{\otimes n}\, & \Delta_n(x) = \sum_{(x)} x_{(1)} \otimes \dots \otimes x_{(n)}
\end{align*}
for iterated comultiplications, where in particular $\Delta_0 = \epsilon$.

The existence of the $\K$-monoid morphism $\Delta$ alone suffices for defining actions on tensor products of $\h$-modules, 
by taking
$x (v \otimes w) = \sum_{(x)} x_{(1)} v \otimes x_{(2)} w$ for $x \in \h$, $v \in V$, and $w \in W$.
Coassociativty then ensures that the $\K$-isomorphism $(U \otimes V) \otimes W \cong U \otimes (V \otimes W)$ is in fact an $\h$-isomorphism.
The counit is responsible for giving a canonical “trivial” action on $\K$, 
or indeed any $\K$-space.
Since $\K$ is necessarily simple, we fix $S_1 = \K$.
Counitality ensures $V \otimes \K \cong_{\h} V \cong_{\h} \K \otimes V$
for any $\h$-module $V$.
In an arbitrary $\h$-module $V$, a vector $v \in V$ such that $xv = \epsilon(x)v$ for any $x \in \h$ is called an \emph{$\h$-invariant},
and the $\h$-submodule of all such vectors is denoted $V^\h$.
The tensor product of modules endows $\mathrm G_0 \h$ with the structure of a ring, with multiplication $[V][W] = V \otimes W$ and unit $[S_1]$.
The dimension map $\dim : \mathrm G_0 \mathcal H \to \mathbb Z$ is then a ring augmentation.

A $\K$-bimonoid $\h$ enjoys an even more well-behaved representation theory when $\h$ is \emph{Hopf}, and thus called a \emph{Hopf $\K$-algebra}.
These are those $\K$-bimonoids for which there exists a (necessarily unique) linear map $\sigma : \h \to \h$ such that for any $x \in \h$
\begin{equation} \label{eq:antipous}
	\sum_{(x)} \sigma(x_{(1)})x_{(2)} = \epsilon(x) = \sum_{(x)} x_{(1)} \sigma(x_{(2)}),
\end{equation}
called an \emph{antipous}.\footnote{From Greek \emph{\textgreek{ἀντίπους}}, the singular of the more familiar antipodes, although many authors use the back-formation \emph{antipode}.}
It can be shown that any such $\sigma$ is in fact a $\K$-bimonoid antiautomorphism, i.e. reverses order of both multiplication and comultiplication~\cite[Proposition~4.0.1]{sweedlerHopfAlgebras1969}.
Antipodes in some sense generalise inverses,
for example if $G$ is a group then the $\K$-bimonoid $\K G$ defined in \cref{ex:monoid} has an antipous defined by $\sigma(x) = x^{-1}$ for $x \in G$.
Similarly, $\K^G$ has an antipous defined by $\sigma(f)(x) = f(x^{-1})$ for $f \in \K^G$.
Unlike in these “classical” examples, the antipous of a Hopf $\K$-algebra need not be an involution,
although it must have finite order~\cite{radfordOrderAntipodeFinite1976}.
The following example is neither commutative, cocommutative, nor semisimple.

\begin{example}[See \cite{chenIrreducibleRepresentationsClass2000,benkartMcKayMatricesFinitedimensional2022}] \label{ex:double-taft}
	Let $\K$ be an algebraically closed field, $n \geq 2$ an integer,
	and $\zeta_n \in \K$ a primitive $n$th root of unity.
	The \emph{Drinfel'd double of the Taft algebra} $\D_n(\zeta_n)$
	is an $n^4$-dimensional Hopf $\K$-algebra
	generated as a $\K$-monoid by $\{a,b,c,d\}$ such that
	\begin{align*}
		ba &= qab,
		&
		db &= qbd,
		&
		ca &= qac,
		&
		dc &= qcd,
		\\
		bc &= cb,
		&
		da-qad &= 1 - bc,
		&
		a^n &= 0 = d^n,
		&
		b^n &= 1 = c^n
	\end{align*}
	with coproduct and counit
	\begin{align*}
		\Delta(a) &= a \otimes b + 1 \otimes a,
		&
		\Delta(d) &= d \otimes c + 1 \otimes d,
		\\
		\Delta(b) &= b \otimes b,
		&
		\Delta(c) &= c \otimes c,
		\\
		\epsilon (a) &= 0 = \epsilon(d),
		&
		\epsilon (b) &= 1 = \epsilon(c)
	\end{align*}
	and antipous of order $2n$
	\begin{align*}
		\sigma(a) &= -ab^{-1},
		&
		\sigma(b) &= b^{-1},
		&
		\sigma(c) &= c^{-1},
		&
		\sigma(d) &= -dc^{-1}.
	\end{align*}
\end{example}

Fixing a $\K$-monoid antiendomorphism $\sigma$ is already sufficient to endow $V^*$ with the structure of an $\h$-module by taking $(x f)(v) = f(\sigma(x) v)$ for any $x \in \mathcal H$, $v \in V$, and $f \in V^*$.
We call the resulting $\h$-module the \emph{left dual}, denoted simply by $V^*$.
Since $\sigma^{-1}$ is also a $\K$-monoid antiendomorphism, $(xf)(v) = f(\sigma^{-1}(x)v)$ also defines an $\h$-module, called \emph{right dual} and denoted $\leftidx^*V$.\footnote{Some texts use left and right in the opposite way. 
Our terminology is chosen to agree with left and right adjoints in a bicategory, of which these are a special case.}
One then has natural isomorphisms 
\begin{align*}
	\leftidx^*(V^*) &\cong_\h V \cong_\h (\leftidx^*V)^*,
	&
	(V \otimes W)^* &\cong_\h W^* \otimes V^*, 	
	&
	\leftidx^*(V \otimes W) \cong_\h \leftidx^*W \otimes \leftidx^* V.
\end{align*}
\Cref{eq:antipous} ensures that the evaluation $V^* \otimes V \to \K$ and coevaluation $\K \to V \otimes V^*$ maps are $\h$-morphisms.
The operations ${?}^*$ and $\leftidx^*{?}$ induce mutually inverse antiautomorphisms on the augmented ring $\mathrm G_0\h$ and group $\mathrm K_0\h$, in each case permuting the basis elements.

\begin{prop}[{\cite[Lemma A.3]{grinbergCriticalGroupsHopf2020}}] \label{thm:hom-isos}
	Let $\h$ be a Hopf $\K$-algebra and $U$, $V$, and $W$ be $\h$-modules.
	Then we have the natural isomorphisms
	\begin{align}
		\lMod{\h}(U \otimes V, W) &\cong_\K \lMod\h(U, W \otimes V^*),\label{eq:hom-isos1}  \\
		\lMod\h(V^* \otimes U, W) &\cong_\K \lMod\h(U, V \otimes W), \label{eq:hom-isos2} \\
		\lMod\h(U \otimes \leftidx^*V, W) &\cong_\K \lMod\h(U, W \otimes V), \label{eq:hom-isos3} \\
		\lMod\h(V \otimes U, W) &\cong_\K \lMod\h(U, \leftidx^*V \otimes W). \label{eq:hom-isos4}
	\end{align}
\end{prop}

For a $\K$-bimonoid $\h$, the algebra homomorphisms $\h \to \K$ form a monoid $\opn{Gr}(\h^*)$ with unit $\epsilon$
under the multiplication $(\chi \eta)(x) = \sum_{(x)} \chi(x_{(1)}) \eta (x_{(2)})$ for $\chi,\eta \in \opn{Gr}(\h^*)$ and $x \in \h$,
which we call the \emph{dual monoid}.
This monoid may be identified with the (isomorphism classes of) 1-dimensional $\h$-modules under the tensor product, where we denote the $\h$-module corresponding to $\chi \in \opn{Gr}(\h^*)$ by $\K_\chi$.
When $\h$ is Hopf, $\opn{Gr}(\h^*)$ becomes the \emph{dual group} where the inverse of $\chi$ is $\chi^* : x \mapsto \chi (\sigma(x))$,
and in particular $\opn{Gr}(\K G^*) = G^\vee$ is the ordinary dual group of linear characters.

\section{Burnside's theorem for $\K$-bimonoids}

For the remainder of the paper, $\h$ is a $\K$-bimonoid and $V$ is an $\h$-module.
In order to characterise the connectedness of the McKay quiver of $\h$ at $V$, we first give an analogue of \cref{thm:burnside-classical} in this generality.
Since we are concerned with the simple $\h$-modules appearing merely as \emph{subquotients} rather than as submodules, 
we arrive at something slightly different to the generalised Burnside's theorem presented in \cite{rieffelBurnsidesTheoremRepresentations1967,passmanBurnsidesTheoremHopf1995}.
For the hypothesis, it is clear that faithfulness in the sense of $\lAnn{\h}(V) = 0$ is far too strict.

\begin{example}
	Let $\K$ be a field with $\opn{char}\K \neq 2$,
	and $\mathrm C_2 = \langle \zeta :\mid: \zeta^2 = 1 \rangle$ be the cyclic group of order $2$.
	Let $V = \K$ with the $\K \mathrm C_2$-action defined by $\zeta v = -v$.
	Then $V$ is faithful on $\mathrm C_2$ but $\zeta + 1 \in \lAnn{\K \mathrm C_2}(V)$.
\end{example}

Since the sum of biideals is a biideal, any subspace $S \leq_\K \h$ of a $\K$-bialgebra contains a unique largest biideal, which we denote $S^\circ$.
Note that if $\h$ is a Hopf $\K$-algebra, all its biideals are automatically Hopf ideals by \cite{nicholsQuotientsHopfAlgebras1978}.
We call $\Ann^\circ_\h (V) := \Ann_\h (V)^\circ$ the \emph{bimonoid annihilator} of $V$.
Variants of this definition are considered in \cite{banicaHopfImagesInner2010,passmanBurnsidesTheoremHopf1995}.
In particular, a module $V$ for which $\lAnn\h^\circ(V) = 0$ is called \emph{inner faithful}.

\begin{lemma} \label{thm:coidealann}
	Suppose $I$ is a coideal annihilating $\h$-modules $V$ and $W$.
	Then $I$ annihilates $V \otimes W$.
\end{lemma}

\begin{proof}
	Let $x \in I$. Then $\Delta(x) = x_1 + x_2$ for some $x_1 \in I \otimes \h$ and $x_2 \in \h \otimes I$,
	and thus $x(V \otimes W) = x_1(V \otimes W) + x_2 (V \otimes W) = 0$. Thus $I$ annihilates $V \otimes W$.
\end{proof}

\begin{prop} \label{thm:biann}
Let $V$ be a finite-dimensional $\h$-module. Then
\begin{align*}
	\lAnn\h^\circ(V) = \bigcap_{n \in \mathbb{N}} \lAnn\h(V^{\otimes n}).
\end{align*}
\end{prop}

\begin{proof}
	Let $I$ denote the right hand side.
	To show that $\Ann_\h^\circ(V) = I$
	we must show that $I$ is a biideal,
	and that any biideal contained in $\Ann_\h(V)$ is also contained in $I$.
	
	First we show the former.
	As the intersection of ideals, $I$ itself is an ideal.
	Since $I \subseteq \Ann_\h(V^{\otimes 0}) = \ker (\epsilon)$, we have $\epsilon(I) = 0$.
	Let $x \in I$.
	Then $0 = x  V^{\otimes 2n} = \sum_{(x)} x_{(1)} V^{\otimes n}  \otimes x_{(2)} V^{\otimes n}$ for all $n \in \mathbb N$ and thus $\Delta(x) \in \lAnn\h(V^{\otimes n}) \otimes \h + \h \otimes \lAnn\h(V^{\otimes n})$ for all $n \in \mathbb N$, so $\Delta(x) \in I \otimes \h + \h \otimes I$.
	Therefore $I$ is a coideal, and hence a biideal.
	
	To show that any biideal $J$ contained in $\lAnn\h(V)$ is also contained in $I$,
	we prove that any such $J$ annihilates $V^{\otimes n}$ for all $n \in \mathbb N$.
	The $n=1$ case is given, the $n = 0$ case follows from $\epsilon(J) = 0$.
	All other cases follow inductively by \cref{thm:coidealann}.
\end{proof}

The quotient $\K$-bimonoid $\h_V := \h / \lAnn\h^\circ(V)$ with projection $\pi : \h \twoheadrightarrow \h_V$ is equipped with an action on $V$
determined by $\pi (x) \, v = x  v$ for $x \in \h$ and $v \in V$,
where it follows $\lAnn{\h_V}^\circ(V) = 0$.
These data are characterised by the following universal property:

\begin{prop}
	Let $\h'$ be a quotient $\K$-bimonoid of $\h$ with projection $\phi : \h \twoheadrightarrow \h'$,
	acting on $V$ so that $\phi(x)v = xv$ for all $x \in \h$ and $v \in V$.
	Then there exists a unique map $\pi_\phi : \h' \twoheadrightarrow \h_V$ such that $\pi = \pi_\phi \phi$.
	\begin{equation*}
		\begin{tikzcd}
		  {\h} && {\h_V} \\
		  \\
		  {\h'} && {\End_\K(V)}
		  \arrow["\pi", two heads, from=1-1, to=1-3]
		  \arrow["\phi"', two heads, from=1-1, to=3-1]
		  \arrow[from=1-3, to=3-3]
		  \arrow["\pi_\phi",dashed, from=3-1, to=1-3]
		  \arrow[from=3-1, to=3-3]
		\end{tikzcd}
	\end{equation*}
\end{prop}

\begin{proof}
	Since $\ker \phi \subseteq \lAnn\h(V)$ is a biideal, it follows that $\ker \phi \subseteq \lAnn\h^\circ(V)= \ker \pi$. The existence of such a factorization follows.
\end{proof}

\begin{thm} \label{thm:burnside-sharp}
	Let $S_i$ be a simple $\h$-module.
	Then there exists an $n \in \mathbb N$ for which $S_i$ is a subquotient of $V^{\otimes n}$ if and only if $\Ann_\h^\circ(V) S_i = 0$.
\end{thm}

\begin{proof}
	We begin by proving the simpler case that if $\lAnn\h^\circ(V) = 0$,
	then for any $S_i$ there exists an $n \in \mathbb N$ such that $[V^{\otimes n} : S_i] \geq 1$.
	To this end, suppose $\lAnn\h^\circ(V) = 0$ and fix a simple $\h$-module $S_i$.
	Then $d_i[V^{\otimes n} : S_i] = \langle P_i, V^{\otimes n} \rangle 
	=\dim \lMod\h(P_i, V^{\otimes n})$,
	so it suffices to exhibit a nonzero $\h$-morphism $f : P_i \to V^{\otimes n}$ for some $n$.
	Fixing a nonzero $v_0 \in P_i$ and an inclusion $\iota : P_i \hookrightarrow \h$ by way of \cref{eq:regular remak}
	gives a nonzero $\tilde v_0 := \iota(v) \in \h$.
	By \cref{thm:biann}, there exists some $n \in \mathbb N$ and $w \in V^{\otimes n}$ with $\tilde v_0w \neq 0$.
	Now let
	\begin{align*}
		f : P_i &\to V^{\otimes n} \\
		v &\mapsto \iota(v)w,
	\end{align*}
	which is an $\h$-morphism since it is clearly $\K$-linear and $f(xv) = \iota(x v)w = x\iota(v)w$ for any $x \in \h$ and $v \in P_i$.
	Moreover it is nonzero since $f(v_0) = \tilde v_0w \neq 0$ by construction.
	Therefore, if $\lAnn\h^\circ(V) = 0$, then for any $S_i$ there exists an $n \in \mathbb N$ such that $[V^{\otimes n}:S_i] \geq 1$.
	
	Now we prove the theorem in full strength.
	For the forward implication, suppose $[V^{\otimes n} :S_i] \geq 1$.
	Then $S_i$ is a subquotient of $V^{\otimes n}$, so $\lAnn\h(V^{\otimes n}) \subseteq \lAnn\h(S_i)$.
	Since $\lAnn\h^\circ(V)$ annihilates $V^{\otimes n}$,
	it also annihilates $S_i$, as claimed.
	
	For the converse, suppose $\lAnn\h^\circ(V)$ annihilates $S_i$.
	Then both $V$ and $S_i$ are naturally regarded as $\h_V$-modules.
	Since $\lAnn{\h_V}^\circ(V) = 0$, it follows from the simpler case that $S_i$ is a subquotient of some $V^{\otimes n}$ as an $\h_V$-module,
	and thus also as an $\h$-module, as claimed.
\end{proof}

Following \cite{grinbergCriticalGroupsHopf2020} we call $V$ \emph{tensor-rich} if for every simple $\h$-module $S_i$ there exists an $n \in \mathbb N$ such that $[V^{\otimes n}: S_i ] \geq 1$.
Recall that the \emph{Jacobson radical} of $\h$ is given by $\Jac(\h) = \bigcap_{i=1}^t \lAnn\h(S_i)$. Therefore:

\begin{cor} \label{thm:burnside-cor}
	$V$ is tensor-rich if and only if $\Ann_\h^\circ(V) \subseteq \Jac(\h)$.
\end{cor}

\section{McKay quivers of $\K$-bimonoids}

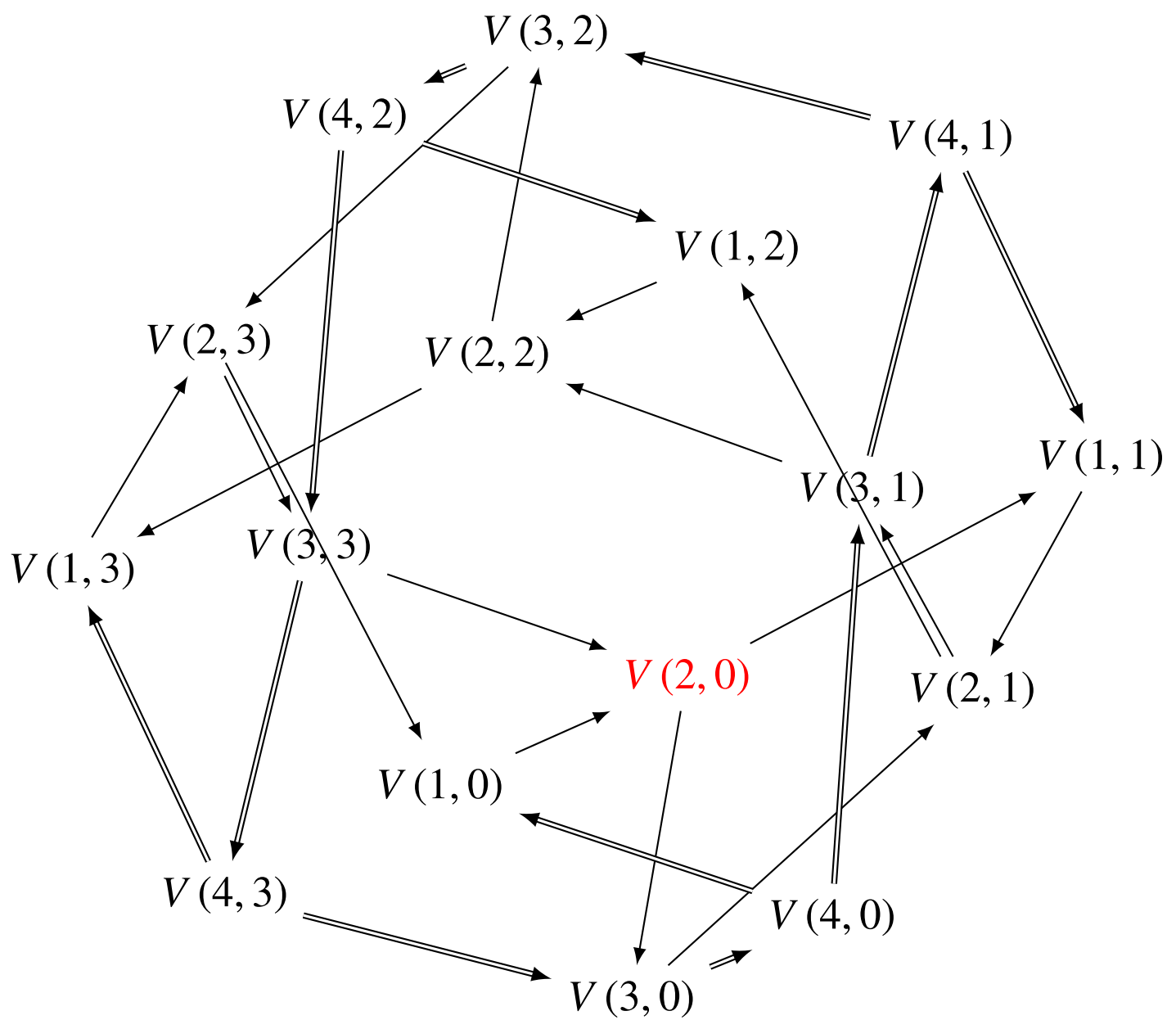
\begin{figure}
	\[
		\begin{tikzpicture}[>=latex,line join=bevel,]
		\node (node_0) at (98.652bp,56.868bp) [draw,draw=none] {$V\left(1, 0\right)$};
		  \node (node_4) at (155.24bp,82.029bp) [draw,draw=none] {$\color{red}V\left(2, 0\right)$};
		  \node (node_1) at (250.03bp,132.29bp) [draw,draw=none] {$V\left(1, 1\right)$};
		  \node (node_5) at (220.55bp,79.022bp) [draw,draw=none] {$V\left(2, 1\right)$};
		  \node (node_2) at (166.55bp,180.2bp) [draw,draw=none] {$V\left(1, 2\right)$};
		  \node (node_6) at (109.33bp,155.86bp) [draw,draw=none] {$V\left(2, 2\right)$};
		  \node (node_3) at (14.39bp,106.32bp) [draw,draw=none] {$V\left(1, 3\right)$};
		  \node (node_7) at (45.632bp,158.78bp) [draw,draw=none] {$V\left(2, 3\right)$};
		  \node (node_8) at (142.51bp,8.3018bp) [draw,draw=none] {$V\left(3, 0\right)$};
		  \node (node_9) at (195.27bp,124.92bp) [draw,draw=none] {$V\left(3, 1\right)$};
		  \node (node_10) at (122.81bp,229.55bp) [draw,draw=none] {$V\left(3, 2\right)$};
		  \node (node_11) at (68.525bp,111.85bp) [draw,draw=none] {$V\left(3, 3\right)$};
		  \node (node_12) at (188.51bp,27.005bp) [draw,draw=none] {$V\left(4, 0\right)$};
		  \node (node_13) at (215.47bp,205.53bp) [draw,draw=none] {$V\left(4, 1\right)$};
		  \node (node_14) at (76.725bp,210.41bp) [draw,draw=none] {$V\left(4, 2\right)$};
		  \node (node_15) at (49.315bp,32.126bp) [draw,draw=none] {$V\left(4, 3\right)$};
		  \draw [black,->] (node_0) -- (node_4);
		  \draw [black,->] (node_1) -- (node_5);
		  \draw [black,->] (node_2) -- (node_6);
		  \draw [black,->] (node_3) -- (node_7);
		  \draw [black,->] (node_4) -- (node_1);
		  \draw [black,->] (node_4) -- (node_8);
		  \draw [black,->,transform canvas={xshift=-3pt}] (node_5) -- (node_2);
		  \draw [black,->] (node_5) -- (node_9);
		  \draw [black,->] (node_6) -- (node_3);
		  \draw [black,->] (node_6) -- (node_10);
		  \draw [black,->,transform canvas={yshift=3pt}] (node_7) -- (node_0);
		  \draw [black,->] (node_7) -- (node_11);
		  \draw [black,->] (node_8) -- (node_5);
		  
		  \draw [black,double,->] (node_8) -- (node_12);
		  
		  \draw [black,->] (node_9) -- (node_6);
		  
		  \draw [black,double,->] (node_9) -- (node_13);
		  
		  \draw [black,->] (node_10) -- (node_7);
		  
		  \draw [black,double,->] (node_10) -- (node_14);
		  
		  \draw [black,->] (node_11) -- (node_4);
		  
		  \draw [black,double,->] (node_11) -- (node_15);
		  
		  \draw [black,double,->] (node_12) -- (node_0);
		  
		  \draw [black,double,->] (node_12) -- (node_9);
		  
		  \draw [black,double,->] (node_13) -- (node_1);
		  
		  \draw [black,double,->] (node_13) -- (node_10);
		  
		  \draw [black,double,->] (node_14) -- (node_2);
		  
		  \draw [black,double,->] (node_14) -- (node_11);
		  
		  \draw [black,double,->] (node_15) -- (node_3);
		  
		  \draw [black,double,->] (node_15) -- (node_8);
		\end{tikzpicture}
	\]
	\caption{
		The McKay quiver of $\mathcal D_4(\zeta_4)$ as in \cref{ex:double-taft} at the module marked in red.
		Simple modules are named according to their classification in \cite{chenIrreducibleRepresentationsClass2000}, where the first index is the dimension.
		The quiver was constructed using the tensor product decompositions given in \cite{chenIndecomposableDecompositionTensor2017}.
	}
	\label{fig:double-taft}
\end{figure}

We define the McKay quiver $\Gamma_\h(V)$ of $\h$ at $V$ by its adjacency matrix $\mathbf M_\h(V) = (M_{ij})$ where
\begin{align*}
	M_{ij} := [S_i \otimes V : S_j]	= \frac{1}{d_j} \langle P_j, S_i V \rangle.
\end{align*}
Taking $\h = \K G$ recovers the McKay quiver of $G$ at $V$, hereinafter called a \emph{classical McKay quiver}.
A non-classical McKay quiver  is shown in \cref{fig:double-taft}.
The following example gives a sense in which classical McKay quivers are dual to Cayley graphs of groups.

\begin{example}\label{ex:cayley}
	Let $M$ be a finite monoid, and $\K^M$ be the $\K$-bimonoid of functions defined in \cref{ex:monoid}.
	For $x \in M$, let $\delta_x \in \K^M$ be the function defined by $\delta_x(y) = 1$ if $y = x$ and $0$ otherwise.
	Since $\K^M = \bigoplus_{x \in M} \K\delta_x$, we see that $\K^M$ is semisimple.
	Each simple module $S_x := \mathbb K \delta_x$ then has the $\K^M$-action $f \delta_x = f(x) \delta_x$.
	Thus all $\K^M$-modules are isomorphic to $V_A := \bigoplus_{x \in A} S_x$
	for some multiset $A \subseteq M$.
	The comultiplication of $\K^M$ induces a $\K^M$-action on $S_x \otimes S_y$ so that $S_x \otimes S_y \cong_{\K^M} S_{x y}$,
	whence we have an isomorphism $\mathrm G_0 \K^M \cong \Z M$ of rings.
	The McKay quiver of $\K^M$ at $V_A$ is precisely the Cayley graph of $M$ with joining set $A$:
	An edge from $S_x$ to $S_y$ for each occurance of $S_y$ in $S_x \otimes V_A$ amounts to an edge from $S_x$ to $S_{xz}$ for each $z \in A$.
\end{example}

\begin{rem}
	More generally, if $\h$ is commutative and $\K$ is a splitting field,\footnote{If in addition $\h$ is semisimple, then these conditions imply $\h \cong \K^M$ as in \cref{ex:cayley} for $M = \opn{Gr}(\h^*)$ by \cite[Theorem 2.3.1]{montgomeryHopfAlgebrasTheir1993}, but this can fail in general.} 
	then all simple modules are 1-dimensional and $\Gamma_\h(V)$ will be the Cayley graph of the dual monoid $\opn{Gr}(\h^*)$,
	where the joining multiset is the Jordan-Hölder decomposition of $V$.
\end{rem}

\begin{prop}
	Let $\h = \K^M$ and consider $V_A$ for $A \subseteq M$ as in \cref{ex:cayley}.
	Then
	\begin{align*}
		\lAnn\h^\circ(V_A) = \opn{span}_\K\{\delta_x : x \in M \setminus \langle A \rangle\}
	\end{align*}
	where $\langle A \rangle$ is the submonoid of $M$ generated by $A$.
\end{prop}
\begin{proof}
	Note that $V_A^{\otimes n} = \left(\bigoplus_{x \in A} S_a\right)^{\otimes n} = \bigoplus_{x \in A^n} S_{\mu_n(x)}$, 
	where $\mu_n : M^n \to M$ denotes the $n$-ary multiplication in $M$.
	Therefore $\lAnn\h\left(V_A^{\otimes n}\right) = \opn{span}_\K\{\delta_x : x \in M \setminus \mu_n(A^n)\}$,
	whence $\lAnn\h^\circ(V_A) = \opn{span}_\K\{\delta_x : x \in M \setminus \langle A \rangle\}$ by \cref{thm:biann}.
\end{proof}

By our discussion in \cref{sec:prelim}, we have an isomorphism $\mathrm G_0 \h \cong_\Z \Z^t$
sending the basis $[S_1], \dots, [S_t]$ to the column vectors which we denote $\mathbf e_1 , \dots, \mathbf e_t$.
We also let
\begin{align*}
	\mathbf s &:= [\dim S_1, \dots, \dim S_t]^\top, & \mathbf p &:= [\dim P_1, \dots, P_t], & \mathbf D &:= \opn{diag}(d_1, \dots, d_t), 
\end{align*}
so that $\langle P_i, S_j \rangle = \mathbf e_i^\top \mathbf D \mathbf e_j$.
If $\mathbf w \in \Z^t$ is a column vector corresponding to $[W] \in \mathrm K_0 \h$,
we see that $\mathbf M_\h(V)^\top \mathbf w$ is the column vector corresponding to $[W][V]$.
Thus $\mathbf M_\h(V)^\top$ is a matrix representation of right multiplication by $[V]$ in $\mathrm K_0 \h$.

\begin{lemma} \label{thm:walks}
	The number of walks of length $n \in \mathbb N$ from $S_i$ to $S_j$ in $\Gamma_\h(V)$ is $[S_i \otimes V^{\otimes n} : S_j]$.
\end{lemma} 
\begin{proof}
	By elementary graph theory, the number of such walks is
	\begin{equation*}
		\mathbf e_i^\top \mathbf M_\h(V)^n \mathbf e_j = \mathbf e_j^\top (\mathbf M_\h(V)^\top)^n \mathbf e_i,
	\end{equation*}
	which by the discussion above is $[S_i \otimes V^{\otimes n} : S_j]$.
\end{proof}

\Cref{thm:burnside-sharp} says there is a walk from $S_1$ to $S_i$ just whenever $\Ann_\h^\circ(V)S_i = 0$,
so in particular such an $S_i$ is weakly connected to $S_1$.
For classical McKay quivers, weakly and strongly connected components are in fact the same~\cite[Proposition~3.1]{browneConnectivityPropertiesMcKay2021}.
For McKay quivers of arbitrary $\K$-bimonoids this is not the case.

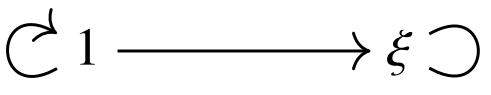
\begin{figure}
	\[\begin{tikzcd}[cramped]
		1 && \xi
		\arrow[from=1-1, to=1-1, loop, in=150, out=210, distance=5mm]
		\arrow[from=1-1, to=1-3]
		\arrow[no head, from=1-3, to=1-3, loop, in=330, out=30, distance=5mm]
	\end{tikzcd}\]
	\caption{\label{fig:cayleyweak} The Cayley graph of the monoid $M = \langle\xi :\mid: \xi^2 = \xi\rangle$ with $M$ as the joining set.}
\end{figure}

\begin{example}
	Let $M$ be the monoid presented by $\langle \xi :\mid: \xi^2 = \xi \rangle$.
	Then $\K^M$ is faithful on the regular module $\K^M$, 
	yet $\Gamma_{\K^M}(\K^M)$, shown in \cref{fig:cayleyweak}, is weakly connected, but not strongly connected.
\end{example}

\Cref{thm:cartan-matrix-symmetry,thm:grinberg-3.6,thm:dimension-eigenvectors} extend results from \cite[\S 3]{grinbergCriticalGroupsHopf2020} to
a possibly non-splitting field $\K$.
\begin{lemma} \label{thm:cartan-matrix-symmetry}
	If $\h$ is Hopf, $d_i[V \otimes P_j^* : S_i] = d_j[\leftidx^*P_i \otimes V : S_j]$ for $i,j \in \{ 1,\dots, t \}$.
	In particular when $V = \K$ we have $d_i[P_j^* : S_i] = d_j[\leftidx^*P_i : S_j]$.
\end{lemma}
\begin{proof}
	From \cref{eq:hom-isos1,eq:hom-isos4} we have
	\begin{equation*}
		d_i[V \otimes P_j^* : S_i] 
		= \langle P_i , V \otimes P_j^* \rangle
		= \langle P_i \otimes P_j, V \rangle
		= \langle P_j, \leftidx^*P_i \otimes V \rangle
		= d_j[\leftidx^*P_i \otimes V : S_j]
	\end{equation*}
	as claimed.
\end{proof}

\begin{lemma} \label{thm:grinberg-3.6}
	If $\h$ is Hopf, for $i = 1,\dots,t$ we have
	\begin{align*}
		[V \otimes P_i^*] = \sum_{j=1}^t	 \frac{d_i}{d_j} [S_j \otimes V : S_i][\leftidx^*P_j] \in \mathrm G_0\h.
	\end{align*}
\end{lemma}

\begin{proof}
	To prove the equality, it suffices to prove the multiplicities of $S_k$ in the left and right elements coincide.
	Applying \cref{thm:cartan-matrix-symmetry} twice
	\begin{align*}
		[V \otimes P_i^* : S_k] 
		&= \frac{d_i}{d_k} [\leftidx^* P_k \otimes V : S_i]
		= \sum_{j=1}^t \frac{d_i}{d_k}[S_j \otimes V : S_i] [\leftidx^*P_k : S_j] \\
		&= \sum_{j=1}^t \frac{d_i}{d_j} [S_j \otimes V : S_i]  [P_j^* : S_k]
	\end{align*}
	as required.
\end{proof}

\begin{lemma} \label{thm:dimension-eigenvectors}
	The vector $\mathbf s$ is a right eigenvector for $\mathbf M_\h(V)$ with eigenvalue $\dim V$.
	In addition, if $\h$ is Hopf,
	then $\mathbf p\mathbf D^{-1}$ is a left eigenvector with eigenvalue $\dim V$.
\end{lemma}

\begin{proof}
	On the one hand,
	\begin{align*}
		(\dim V)\mathbf s 
		&= \sum_{i=1}^t(\dim V)(\dim S_i) \mathbf e_i
		= \sum_{i=1}^t \dim (S_i \otimes V) \mathbf e_i\\
		&= \sum_{i,j=1}^t [S_i \otimes V : S_j](\dim S_j) \mathbf e_i
		= \mathbf M_\h(V) \mathbf s
	\end{align*}
	If $\h$ is Hopf, applying \cref{thm:grinberg-3.6} gives
	\begin{align*}
		(\dim V) \mathbf p	\mathbf D^{-1}
		&= \sum_{i=1}^t \frac{1}{d_i} (\dim V)(\dim P_i) \mathbf e_i^\top 
		= \sum_{i=1}^t \frac{1}{d_i} (\dim V)(\dim P_i^*) \mathbf e_i^\top \\
		&= \sum_{i=1}^t \frac{1}{d_i} \dim (V \otimes P_i^*) \mathbf e_i^\top 
		= \sum_{i,j=1}^t \frac{1}{d_j} [S_j \otimes V : S_i] (\dim P_j) \mathbf e_i^\top \\
		&= \mathbf p \mathbf{D}^{-1} \mathbf M_\h(V)
	\end{align*}
	as claimed.
\end{proof}

\begin{lemma} \label{thm:strong-weak}
	Let $\Gamma$ be a quiver and $\mathbf M$ be its adjacency matrix.
	If $\mathbf M$ has positive left and right eigenvectors of eigenvalue $\lambda$, then the strongly and weakly connected components of $\Gamma$ coincide.
	Moreover, $\lambda$ is the spectral radius of $\mathbf M$ 
	and the unique eigenvalue possessing a positive eigenvector,
	and its algebraic and geometric multiplicity equal the number of connected components of $\Gamma$.
\end{lemma}
\begin{proof}
	Suppose $\mathbf M = (M_{ij})$ has a positive right eigenvector $\mathbf v = \sum_{i \in \Gamma_0} v_i \mathbf e_i$
	and a positive left eigenvector $\mathbf f = \sum_{i \in \Gamma_0} f_i\mathbf e_i^\top$, both of eigenvalue $\lambda$.
	We argue that a set $S \subseteq \Gamma_0$ of vertices has an incoming edge if and only if is has an outgoing edge.
	Letting $S' := \Gamma_0 \setminus S$, on the one hand
	\begin{align*}
		\lambda \sum_{i\in S}f_iv_i = \sum_{i \in S, j \in \Gamma_0} f_i M_{ij}v_j = \sum_{i,j \in S} f_i M_{ij}v_j + \sum_{i\in S, j \in S'} f_iM_{ij}v_j
	\end{align*}
	and on the other
	\begin{align*}
		\lambda	\sum_{j\in S}f_jv_j = \sum_{i \in \Gamma_0, j \in S} f_iM_{ij}v_j = \sum_{i,j \in S} f_i M_{ij}v_j + \sum_{i \in S' ,j \in S} f_iM_{ij} v_j.
	\end{align*}
	Since the left hand sides of these equations are equal, so too are the right hand sides, and thus
	\begin{align} \label{eq:SS'}
		\sum_{i \in S, j \in S'}f_i M_{ij} v_j = \sum_{i \in S', j \in S} f_i M_{ij} v_j.
	\end{align}
	Since $\mathbf f$ and $\mathbf v$ are positive and $\mathbf M$ is nonnegative,
	the left hand side of \cref{eq:SS'} is positive if and only if $M_{ij} > 0$ for some $i \in S$ and $j \in S'$, i.e. $S$ has an outgoing edge,
	while the right hand side is positive if and only if $M_{ij} > 0$ for some $i \in S'$ and $j \in S$, i.e. $S$ has an incoming edge.
	Thus $S$ has an incoming edge if and only it has an outgoing edge.
	
	Now fix some vertex $i \in \Gamma_0$, let $S$ be the set of all $j \in \Gamma_0$ for which there exists a directed path $i \rightsquigarrow j$, and $T$ be the set of all $j \in \Gamma_0$ for which there exists a directed path $j \rightsquigarrow i$.
	Since $S$ has no outgoing edges, it has no incoming edges and thus $T \subseteq S$.
	Similarly, since $T$ has no incoming edges, $S \subseteq T$, 
	so in fact $S = T$ is the strongly connected component of $i$ and has no incoming or outgoing edges.
	Therefore the strongly and weakly connected components of any $i \in \Gamma_0$ coincide.
	
	It follows that up to permutation, $\mathbf M = \bigoplus_{i=1}^n\mathbf M_i$ is block diagonal,
	where $\mathbf M_1, \dots, \mathbf M_n$ are the adjacency matrices of the connected components of $\Gamma$.
	The right eigenvector $\mathbf v = \bigoplus_{i=1}^n\mathbf v_i$ is the direct sum of Perron-Frobenius right eigenvectors for each block.
	Since $\lambda$ is both the spectral radius and an eigenvalue of each block
	with multiplicity 1, it follows that $\lambda$ is the spectral radius and an eigenvalue of $\mathbf M$ with multiplicity $n$.
\end{proof}

\begin{thm} \label{thm:mckay-strong-weak}
	If $\h$ is Hopf, the strongly and weakly connected components of $\Gamma_\h(V)$ coincide.
	Moreover, $\dim V$ is the spectral radius of $\mathbf M_\h(V)$ and the unique eigenvalue possessing a positive eigenvector,
	and its algebraic and geometric multiplicity equal the number of connected components of $\Gamma_\h(V)$.
\end{thm}
\begin{proof}
	Immediate from \cref{thm:dimension-eigenvectors,thm:strong-weak}.
\end{proof}

In light of \cref{thm:mckay-strong-weak} we drop the qualifiers of strongly and weakly when discussing connectedness of McKay quivers of Hopf $\K$-algebras.

\begin{cor} \label{thm:criterion}
	If $\h$ is Hopf, $\Gamma_\h(V)$ is connected if and only if $\lAnn\h^\circ(V) \subseteq \Jac(\h)$.
\end{cor}

The following is a straightforward generalisation of \cite[Lemma 4.1]{browneConnectivityPropertiesMcKay2021} to the Hopf setting.

\begin{prop}
	If $\h$ is Hopf, the dual group $\opn{Gr}(\h^*)$ acts on the vertices of $\Gamma_\h(V)$ as quiver automorphisms,\footnote{We take a quiver automorphism to be a bijection $\phi : \Gamma_0 \to \Gamma_0$ on vertices such that $\Gamma(v,w) \cong \Gamma(\phi(v),\phi(w))$ as sets,
		i.e. the identities of individual edges does not matter.
	}
	where $\chi \in \opn{Gr}(\h^*)$ takes $S$ to $\mathbb K_\chi \otimes S$.
	This action is regular on the vertices corresponding to $1$-dimensional $\h$-modules.
\end{prop}
\begin{proof}
	For any $\h$-module $W$ and $\chi \in \opn{Gr}(\h^*)$, the $\h$-module $\K_\chi \otimes W$ may be viewed as having the same underlying $\K$-space as $W$ 
	but the action $\sum_{(x)} \chi(x_{(1)})\,x_{(2)}v$ for $x \in \h$ and $v \in W$.
	It follows that for $U \leq_\K W$, if $U \leq_\h W$ then also $U \leq_\h \K_\chi \otimes W$.
	Since $\chi^* \chi = \epsilon$, replacing $\chi$ with $\chi^*$ and $W$ with $\K_\chi \otimes W$ in the preceding argument shows the converse, 
	so in fact the $\h$-submodules of $W$ and $\K\chi \otimes W$ coincide.
 	Therefore $W$ is simple if and only if $\chi W$ is simple.
 	
 	It follows that left multiplication by $\chi$ gives a group isomorphism of $\mathrm G_0 \h$ permuting the basis $[S_1],\dots, [S_t]$, whence $[\chi S_i \otimes V : \chi S_j] = [S_i \otimes V : S_j]$,
 	so this action on vertices is indeed a quiver isomorphism.
 	The action of $\opn{Gr}(\h^*)$ on the 1-dimensional $\h$-modules is equivalent to its action on itself by left-multiplication,
 	and is thus regular.
 \end{proof}
 
 \begin{example}
 	The action of the dual group $\opn{Gr}(\mathcal D_4(\zeta_4)^*) \cong \mathrm C_4$ on the quiver in \cref{fig:double-taft} gives all of its automorphisms.
 \end{example}
 
 \begin{cor}
 	If $\h$ is Hopf, the connected components of $\Gamma_\h(V)$ containing 1-di\-men\-sion\-al $\h$-modules are pairwise isomorphic.
 \end{cor}

\section{McKay quivers of groups over modular fields}

To conclude, we apply our results to the case where $\h = \K G$ for $G$ a finite group and $\K$ a splitting field for $G$ and its subgroups with characteristic $p$ dividing $|G|$.
The Jacobson radical of $\K G$ is known, being given by the annihilator of the subspace of $\K G$ spanned by $p$-regular section sums \cite{tsushimaNotesRadicalFinite1978}.
However, the resulting specialisation of \cref{thm:criterion} is more directly obtained by determining the number of connected components exactly using Brauer characters.
Let $\beta_1,\dots,\beta_t$ denote the Brauer characters of $S_i, \dots, S_t$, and $c_i^G, \dots, c_t^G$ enumerate $p$-regular conjugacy classes.
Furthermore let $\beta_V$ denote the Brauer character of $V$ and $\rho_V : G \to \GL(V)$ its induced homomorphism.

\begin{prop} \label{thm:brauer-mckay-eigenvalues}
	The eigenvalues of $\mathbf M_{\K G}(V) = M_{ij}$ are $\beta_V(c_k)$ with corresponding right eigenvectors $[\beta_1(c_k), \dots, \beta_t(c_k)]^\top$ for $k = 1,\dots,t$.
\end{prop}
\begin{proof}
	For each $i=1,\dots,t$ we have $\sum_{j=1}^t  M_{ij} \beta_j = \beta_i \beta_V$ by \cite[Proposition 10.1.3]{webbBrauerCharacters2016}.
	Evaluating at $c_k$ then gives $\sum_{j=1}^t M_{ij} \beta_j(c_k) = \beta_i(c_k)\beta_V(c_k)$,
	so $[\beta_1(c_i), \dots, \beta_t(c_i)]^\top$ is a right eigenvector of eigenvalue $\beta_V(c_k)$.
	By \cite[Proposition 10.2.5]{webbBrauerCharacters2016} these eigenvectors are linearly independent and thus form a complete set.
\end{proof}

\begin{prop} \label{thm:modular-mckay-components}
	The number of connected components in $\Gamma_{\K G}(V)$ equals the number of $p$-regular conjugacy classes contained in $\ker \rho_V$.
\end{prop}
\begin{proof}
	For $p$-regular $g \in G$, we have $\beta_V(g) = \beta(1) = \dim V$ if and only if $g \in \ker \rho_V$ by the proof of \cite[Theorem 1.1]{lewisKernelsBrauerCharacters2026}.
	Therefore, the multiplicity of $\dim V$ as an eigenvalue of $\mathbf M_{\K G}(V)$ equals the number of $p$-regular conjugacy classes contained in $\ker \rho_V$ by \cref{thm:brauer-mckay-eigenvalues},
	which in turn equals the number of connected components in $\Gamma_{\K G}(V)$ by \cref{thm:mckay-strong-weak}.
\end{proof}

\begin{cor}
	$\Gamma_{\K G}(V)$ is connected if and only if $\ker \rho_V$ is a $p$-group.
\end{cor}
\begin{proof}
	If $\ker \rho_V$ is a $p$-group, then $\{1\}$ is only $p$-regular conjugacy class it contains.
\end{proof}

\acks
This paper is based on an Honours project supervised by John Bamberg at the University of Western Australia in 2025.
The author is extremely grateful to John for his feedback and support continuing beyond graduation.

The McKay quivers in \cref{fig:PSL32,fig:double-taft} were generated using SageMath interfacing with GAP.
The relevant Jupyter notebook which includes a number of additional examples is available at \href{https://github.com/jajaperson/sage-mckay-quivers/}{\color{blue}https://github.com/jajaperson/sage-mckay-quivers/}.

\bibliographystyle{srtnumbered}
\bibliography{McKayHopf}

\end{document}